\documentclass [11pt,oneside]{article}

\reversemarginpar \textwidth=14cm \textheight=19cm

\usepackage{amssymb} \usepackage{amsfonts} \usepackage{amsmath}
\usepackage{amsthm} \usepackage{epsfig}
\usepackage{graphics}

\newtheorem{lemma}{Lemma}[section] 
\newtheorem{teo}[lemma]{Theorem}
\newtheorem{rem}[lemma]{Remark} 
\newtheorem{prop}[lemma]{Proposition}
\newtheorem{cor}[lemma]{Corollary} 

\newtheorem{ex}[lemma]{Example} 
\newtheorem{defn}[lemma]{Definition}

\newcommand{\matN}{\ensuremath {\mathbb{N}}}
\newcommand{\matR} {\ensuremath {\mathbb{R}}}

\newcommand{\matZ} {\ensuremath {\mathbb{Z}}}

\newcommand{\matH} {\ensuremath {\mathbb{H}}}

\newcommand{\matE} {\ensuremath {\mathbb{E}}}

\newcommand{\calA} {\ensuremath {\mathcal{A}}}

\newcommand{\calT} {\ensuremath {\mathcal{T}}}

\newcommand{\nota} [1] {\caption{\footnotesize{#1}}}

\newfont{\Got}{eufm10 scaled 1200}

\font\titsc=cmcsc10 scaled 1200

\newcommand{\finedimo}{{\hfill\hbox{$\square$}\vspace{2pt}}}

\author{Roberto \titsc{Frigerio}\footnote{This 
research was partially supported by the University of Melbourne.}}

\title{Hyperbolic manifolds with geodesic boundary\\ which are determined by their fundamental group}

\begin{document}

\maketitle

\begin{abstract}
Let $M_1$ and $M_2$ be $n$-dimensional connected 
orientable finite-volume 
hyperbolic manifolds with geodesic boundary, and let
$\varphi:\pi_1(M_1)\to\pi_1(M_2)$ be a given group
isomorphism.
We study the problem whether there exists an isometry
$\psi:M_1 \to M_2$ such that $\psi_{\ast}=\varphi$.
We show that this is always the case if $n\geqslant 4$, 
while in the $3$-dimensional case the existence of $\psi$ is  
proved under some (necessary) additional conditions on $\varphi$.
Such conditions are trivially satisfied if  
$\partial M_1$ and $\partial M_2$ are both compact.	

\vspace{4pt}

\noindent MSC (2000): 30F40 (primary), 57N16 (secondary).
\end{abstract}

\noindent
Let $M_1$ and $M_2$ be connected orientable finite-volume
hyperbolic $n$-manifolds with geodesic boundary. Suppose 
$n\geqslant 3$ and let $\varphi:\pi_1(M_1)\to\pi_1(M_2)$
be an isomorphism of abstract groups.
%In this paper 
We determine necessary and
sufficient conditions  
for $\varphi$  to be induced by an isometry $\psi:M_1 \to M_2$.
When this is the case, 
we say that $\varphi$ is \emph{geometric}
(see Section~\ref{first:sec} for a more detailed definition).
Mostow-Prasad's rigidity theorem 
ensures geometricity of $\varphi$ whenever  
the boundary of $M_i$ is empty for $i=1,2$. 
Building on a result of Floyd~\cite{Floyd},
we will extend Mostow-Prasad's result to the
non-empty boundary case, following slightly different strategies 
according to the dimension of the manifolds
involved. 

If $M_1$ and $M_2$ are $3$-dimensional hyperbolic
manifolds with non-empty geodesic boundary, 
applying Mostow-Prasad's rigidity theorem to their doubles,
\emph{i.e.} to the manifolds obtained by mirroring
$M_1$ and $M_2$ in their boundary, we will show 
that $\varphi$ is geometric  provided it is induced
by a homeomorphism, rather than an isometry. 
A result of Marden and 
Maskit~\cite{MM} 
will then be applied to relate
the existence of a homeomorphism
inducing $\varphi$ to the behaviour of $\varphi$ with respect
to the \emph{peripheral} subgroups of $\pi_1(M_1)$ and $\pi_1(M_2)$ 
(see below for a definition). 

If $\mathrm{dim}(M_1)=\mathrm{dim}(M_2)\geqslant 4$, 
the existence of an isometry
$\psi:M_1 \to M_2$ such that $\psi_{\ast}=\varphi$ will be proved
by a more direct argument using results from~\cite{Tukia2}.

\section{Preliminaries and statement}\label{first:sec}

In this section we
list some preliminary facts about the topology
and geometry of orientable finite-volume hyperbolic $n$-manifolds
with geodesic boundary and we state our main theorem and its
corollaries. From now on we will always suppose $n\geqslant 3$.
Moreover, all manifolds will be connected and orientable.
We omit all proofs about the basic material addressing the
reader to~\cite{FriPe,Kojima1,Kojima2}.

Before going into the real matter, we
devote the first paragraph to give a formal definition of the notion
of \emph{geometric isomorphism} between fundamental
groups of hyperbolic manifolds. To this aim we will need 
to spell out in detail some well-known elementary results in 
the theory of fundamental groups. 

\paragraph{Homomorphisms between fundamental groups}
If $\varphi,\varphi':G\to H$ are group homomorphisms,
we say that $\varphi'$ is conjugated to
$\varphi$ if there exists $h\in H$ such that
$\varphi'(g)=h\varphi(g)h^{-1}$ for every $g\in G$.
Let $X$ be a manifold and $x_0,x_1$ be points in $X$. Then
there exists an isomorphism $\pi_1(X,x_0)\cong
\pi_1(X,x_1)$ which is canonical up to conjugacy.
It follows that an abstract group $\pi_1(X)$ 
is well-defined and for any $x_0\in X$ there exists
a preferred conjugacy class of isomorphisms
between $\pi_1(X)$ and $\pi_1(X,x_0)$.
%Now let $\pi_X:\widetilde{X}\to X$ be a universal covering of
%$X$ and $\mathrm{Deck}(\pi_X)$ be the group of covering automorphisms of 
%$\pi_X$. 
%Any couple of points $x_0\in X,\widetilde{x}_0\in\pi_X^{-1}
%(x_0)$ determines a natural identification $\pi_1(X,x_0)\cong
%\mathrm{Deck}(\pi_X)$, which in turn gives a natural 
%conjugacy class of isomorphisms between $\pi_1(X)$
%and $\mathrm{Deck}(\pi_X)$.

If $f:X\to Y$ is a continuous map between manifolds,
then $f$ determines a well-defined conjugacy class of homomorphisms
$f_\ast\in\mathrm{Hom}(\pi_1(X),\pi_1(Y))/\pi_1(Y)$.
If a homomorphism $\varphi:\pi_1(X)\to\pi_1(Y)$ is given, we 
say that $\varphi$ is induced by $f$ if $\varphi$ belongs
to $f_\ast$; if so, with an abuse we will write 
$\varphi=f_\ast$, rather than $[\varphi]= f_\ast$.

\begin{defn}
Let $M_1$ and $M_2$ by hyperbolic manifolds with geodesic boundary and
$\varphi:\pi_1(M_1)\to\pi_1(M_2)$ be a group isomorphism.
Then $\varphi$ is \emph{geometric} if $\varphi=\psi_\ast$
for some isometry $\psi:M_1\to M_2$.
\end{defn}

%Now let $\pi_X:\widetilde{X}\to X$ and $\pi_Y:\widetilde{Y}\to Y$
%\appunto{Da qui a fine par si pu\`o togliere}
%be universal coverings of $X$ and $Y$ respectively.
%We denote by $\mathrm{Deck}(\pi_X)$ and $\mathrm{Deck}(\pi_Y)$
%the group of covering automorphisms of
%$\pi_X$ and $\pi_Y$. Recall that there are isomorphisms
%$\pi_1(X)\cong\mathrm{Deck}(\pi_X)$, $\pi_1(Y)\cong
%\mathrm{Deck}(Y)$, which are canonical up to conjugation.
%Any continuous map $f:X\to Y$
%admits a continuous lift $\widetilde{f}:\widetilde{X}\to
%\widetilde{Y}$ such that $\pi_Y\circ\widetilde{f}=f\circ \pi_X$.
%For any $\gamma\in\mathrm{Deck}(\pi_X)$ there exists a
%unique $f_{\#}(\gamma)\in\mathrm{Deck}(\pi_Y)$ 
%such that $\widetilde{f}\circ\gamma=
%f_{\#}(\gamma)\circ\widetilde{f}$. The map $f_\#:
%\mathrm{Deck}(\pi_X)\to \mathrm
%{Deck}(\pi_Y)$ is in fact a group homomorphism,
%and its conjugacy class, which does not depend on the choice of
%$\widetilde{f}$, is actually equal to $f_\ast$ under the
%natural identifications $\pi_1(X)\cong\mathrm{Deck}(\pi_X)$,
%$\pi_1(Y)\cong\mathrm{Deck}(\pi_Y)$ mentioned above.

\paragraph{Natural compactification of hyperbolic manifolds} Let $N$ 
be a complete finite-volume 
hyperbolic $n$-manifold with (possibly empty) 
geodesic boundary (from
now on we will summarize all this information saying just that 
$N$ is hyperbolic). Then $\partial N$, endowed
with the Riemannian metric it inherits from $N$, is a
hyperbolic $(n-1)$-manifold without boundary
(completeness of $\partial N$ is obvious, and the 
volume of $\partial N$ is proved to be finite in~\cite{Kojima1}).
It is well-known~\cite{Kojima1, Kojima2} that $N$ consists of a compact portion
together with some cusps based on Euclidean $(n-1)$-manifolds. 
More precisely, the $\varepsilon$-thin part of $N$
(see~\cite{thu}) consists of cusps of the form $T\times [0,\infty)$,
where $T$ is a compact Euclidean $(n-1)$-manifold 
with (possibly empty) geodesic boundary such that 
$\left(T\times[0,\infty)\right)\cap\partial N=
\partial T\times[0,\infty)$.
A cusp based on a closed Euclidean $(n-1)$-manifold
is therefore disjoint from $\partial N$ and
is called \emph{internal}, while a cusp
based on a Euclidean $(n-1)$-manifold with non-empty
boundary intersects $\partial N$ in one or two
internal cusps of $\partial N$, 
and is called a \emph{boundary cusp}.  
This
description of the ends of $N$ easily implies that $N$ admits
a natural compactification $\overline{N}$ obtained by adding
a closed Euclidean $(n-1)$-manifold for each internal cusp
and a compact Euclidean $(n-1)$-manifold with non-empty 
geodesic boundary
for each boundary cusp. 

When $n=3$, $\overline{N}$ is obtained by adding to $N$
some tori and some closed annuli.
In this case we denote by $\calA_N$ the family
of added closed annuli, and we observe that
no annulus in $\calA_N$ lies on a torus in
$\partial\overline{N}$. Note also that 
$\calA_N=\emptyset$ if $\partial N$ is compact.
A loop $\gamma\in\pi_1(N)$ will be called an
\emph{annular cusp loop} if it is freely-homotopic 
to a loop in some annulus of $\calA_N$.

\paragraph{Main result} 
We are now ready to state our main result.

\begin{teo}\label{main:teo}
Let $N_1$ and $N_2$ be hyperbolic $n$-manifolds, and let $\varphi:
\pi_1(N_1)\to\pi_1(N_2)$ be a group isomorphism.
If $n=3$, suppose also that the 
following condition holds:
\begin{itemize}
\item
$\varphi(\gamma)$ is an annular cusp loop in $\pi_1(N_2)$
if and only if $\gamma$ is an annular cusp loop in
$\pi_1(N_1)$.
\end{itemize}
Then $\varphi$ is geometric.
\end{teo}

Theorem~\ref{main:teo} readily implies the following corollaries:

\begin{cor}\label{cor1}
Let $N_1$ and $N_2$ be hyperbolic $3$-manifolds with compact geodesic
boundary and let $\varphi:\pi_1(N_1)\to
\pi_1(N_2)$ be an isomorphism. Then $\varphi$
is geometric.
\end{cor}

\begin{cor}\label{cor2}
Let $N$ be a hyperbolic $n$-manifold,
let 
$\mathrm{Iso}(N)$ be the group of isometries of $N$ and
let ${\mathrm{Out}(\pi_1(N)):=}{\mathrm{Aut}(\pi_1(N))/\pi_1(N)}$ 
be the group of the outer
isomorphisms of $\pi_1(N)$.
If $n=3$, suppose also that the boundary of $N$ is compact.
Then there is a natural isomorphism $\mathrm{Iso}(N)\cong
\mathrm{Out}(\pi_1(N))$.
\end{cor}
\noindent\emph{Proof:}
Let $h:\mathrm{Iso}(N)\to \mathrm{Out}(\pi_1(N))$ be the map
defined by $h(\psi)=\psi_\ast$. Then $h$ is a well-defined homomorphism.
Injectivity of $h$ is a well-known fact, while surjectivity of $h$
is an immediate consequence of Theorem~\ref{main:teo}
and Corollary~\ref{cor1}.
\finedimo

%\paragraph{Hyperbolization}
%Let $\overline{M}$ be a compact orientable manifold 
%with non-empty boundary, let $\calT$ be the set of
%boundary tori of $\overline{M}$ and
%let $\calA$ be a family of disjoint annuli in 
%$\partial \overline{M}\setminus\calT$. 
%Thurston's hyperbolization theorem for Haken manifolds~\cite{thu2}
%gives necessary and sufficient topological conditions
%on the pair $(\overline{M},\calA)$
%for $M=\overline{M}\setminus(\calT\cup\calA)$ to be hyperbolic:

%\begin{teo}\label{hyper:teo}
%The manifold $M=\overline{M}\setminus (\calT\cup\calA)$ is hyperbolic
%if and only if the following conditions hold:
%\begin{itemize}
%\item
%the components of $\partial M$ have negative Euler characteristic;
%\item
%$\overline{M}$ is boundary-irreducible and geometrically atoroidal;
%\item
%the only proper essential annuli contained in $M$ are parallel
%in $\overline{M}$ to the annuli in $\calA$.
%\end{itemize}
%\end{teo}

\paragraph{Universal covering and action at the infinity}
Let $N$ be a $n$-dimensional hyperbolic manifold
and let $\pi:\widetilde{N}\to N$ be
the universal covering of $N$. By developing $\widetilde{N}$
in $\matH^n$ we can identify $\widetilde{N}$ with
a convex polyhedron of $\matH^n$ bounded by
a countable number of disjoint geodesic hyperplanes $S_i,\ i\in\matN$.
For any $i\in\matN$ let $S_i^+$ denote the closed 
half-space of $\matH^n$ 
bounded by $S_i$ and containing $\widetilde{N}$, let
$S_i^-$ be the closed half-space of $\matH^n$ opposite to
$S_i^+$ and let $\Delta_i$ be the internal part of the closure
at infinity of $S_i^-$. Of course we have $\widetilde{N}=
\bigcap_{i\in\matN} S_i^+$, so denoting by 
$\widetilde{N}_{\infty}$ the closure at infinity of $\widetilde{N}$
we obtain $\widetilde{N}_{\infty}=\partial\matH^n\setminus
\bigcup_{i\in\matN} \Delta_i$. 

The group of the automorphisms of the covering
$\pi:\widetilde{N}\to N$ can be identified in a natural way
with 
a discrete torsion-free subgroup $\Gamma$ of 
$\mathrm{Iso}^+(\matH^n)$ such that
$\gamma(\widetilde{N})=\widetilde{N}$ for any $\gamma\in\Gamma$
and $N\cong\widetilde{N}/\Gamma$.
Also recall that there exists an isomorphism
$\pi_1(N)\cong\Gamma$, which is canonical up to
conjugacy. 
Let $\Lambda(\Gamma)$ denote the limit set of $\Gamma$
and let $\Omega(\Gamma)=\partial \matH^n\setminus\Lambda(\Gamma)$.
Kojima has shown in~\cite{Kojima1} that $\Lambda(\Gamma)=
\widetilde{N}_{\infty}$, so 
the round balls $\Delta_i,\ i\in\matN$ previously
defined actually are the connected components of $\Omega(\Gamma)$.
A subgroup of $\Gamma$ is called \emph{peripheral} if it is equal
to the stabilizer of one of the $\Delta_i$'s. 

Since $\widetilde{N}_{\infty}=\Lambda(\Gamma)$,
we have that $\widetilde{N}$
is the intersection of $\matH^n$ 
with the convex hull of $\Lambda(\Gamma)$, so
$N$ is 
the convex core (see~\cite{thu}) of the hyperbolic manifold
$\matH^n/\Gamma$. This implies that $N$ uniquely determines
$\Gamma$ up to conjugation by elements in
$\mathrm{Iso}^+(\matH^n)$, that $\Gamma$ is geometrically finite
and that $N$ is homeomorphic to the manifold
$\left(\matH^3\cup\Omega(\Gamma)\right)/\Gamma$.

\paragraph{Parabolic subgroups of $\Gamma$}
Let $\Gamma'$ be a subgroup of
$\Gamma$. We say that $\Gamma'$ is \emph{maximal parabolic} if it is
parabolic (\emph{i.e.} all its non-trivial elements are parabolic)
and it is maximal with respect to inclusion among 
parabolic subgroups of $\Gamma$.
If $\Gamma'$ is a maximal parabolic subgroup of $\Gamma$, then there
exists a point $q\in\partial\matH^n$ such that $\Gamma'$ equals the 
stabilizer of $q$ in $\Gamma$. Then $\Gamma'$ can be naturally identified
with a discrete subgroup of $\mathrm{Iso}^+(\matE^{n-1})$, so by
Bierbebach's Theorem~\cite{rat} $\Gamma'$ contains an Abelian subgroup $H$ of finite
index. If $k$ is the rank of $H$, we say that $\Gamma'$ is a rank-$k$
parabolic subgroup of $\Gamma$.
Now it is shown in~\cite{Kojima1} that if $i\neq j$, then
$\overline{\Delta}_i\cap\overline{\Delta}_j$ is either empty
or consists of one point $p$ whose stabilizer is a rank-$(n-2)$
parabolic subgroup of $\Gamma$. 
Moreover, any maximal rank-$(n-2)$ parabolic subgroup of $\Gamma$
is the stabilizer of a point $p$ which lies on the boundary
of two different $\Delta_i$'s.
On the other hand,
the intersection of $\widetilde{N}$ with a horoball centered at a 
point with rank-$(n-2)$ parabolic stabilizer
projects onto a boundary cusp of $N$, and
any boundary cusp of $N$ lifts to the intersection
of $\widetilde{N}$ with a horoball centered at a
point with rank-$(n-2)$ parabolic stabilizer.
It follows that there is a natural correspondence
between the boundary cusps of $N$ and the 
conjugacy classes of rank-$(n-2)$ maximal parabolic
subgroups of $\Gamma$. 

We shall see that
rank-$1$ maximal parabolic subgroups of $\Gamma$
play a special role 
in the proof of our main theorem.
Since any parabolic subgroup 
of $\Gamma$ corresponds to a cusp of $N$, we have that
if $n\geqslant 4$
then $\Gamma$ does not contain rank-$1$ maximal parabolic subgroups,
while when $n=3$ the elements of  
rank-$1$ maximal parabolic subgroups of $\Gamma$ 
correspond to
the annular cusp loops previously defined.
For later purpose we point out the following:

%If $C\subset N$ is a toric cusp, then there exist a point $p\in
%\widetilde{N}_{\infty}$ and a horoball $O_p$ centered at 
%$p$ such that $C=\pi(O_p)=O_p/\Gamma_p$, where 
%$\Gamma_p=\{\gamma\in\Gamma :\ \gamma(p)=p\}$ is a
%rank-2 parabolic subgroup of $\Gamma$.
%Moreover, every rank-2 parabolic subgroup
%of $\Gamma$ arises in this way, so toric cusps of $N$
%correspond in a natural way to conjugacy classes
%of rank-2 parabolic subgroups of $\Gamma$.

\begin{rem}\label{nomaximal:rem}
\emph{For any $k\in\matN$ let $H_k$ be the
stabilizer of $\Delta_k$ in $\Gamma$.
If $i\neq j$, then either 
$\overline{\Delta}_i\cap\overline{\Delta}_j=\emptyset$ and
$H_i\cap H_j=\emptyset$, or 
$\overline{\Delta}_i\cap
\overline{\Delta}_j=\{p\}$
and $H_i\cap
H_j$ is the  
rank-$(n-2)$ parabolic stabilizer of $p$ in $\Gamma$.} 
\end{rem}

\section{Some preliminary lemmas}

%Let $N_1$ and $N_2$ be hyperbolic with non empty boundary, let
%$\pi_i:\matH^3\supset\widetilde{N}_i\to N_i$ be the universal
%coverinf of $N_i$ and let
%$\Gamma_i$ be the (conjugacy class of the) 
%subgroup of $PSL(2,\matC)$ such that 
%$N_i\cong\widetilde{N}_i/\Gamma_i$.
The following result is a slight generalization of Lemma 5.1
in~\cite{KMS}, which is due to J.P.~Otal.
Notations are kept from the preceding
section.

\begin{lemma}\label{connect:lemma}
Let $j:S^{n-2}\to\Lambda(\Gamma)$ be a topological
embedding. Then $\Lambda(\Gamma)\setminus j(S^{n-2})$
is path connected if and only if $j(S^{n-2})=\partial \Delta_l$ for
some $l\in\matN$.  
\end{lemma}
\noindent\emph{Proof:}
Suppose that  $j(S^{n-2})=\partial \Delta_0$. 
Using the upper half-space model of hyperbolic space,
we identify $\partial \matH^n$ with
$(\matR^{n-1}\times\{0\})\cup\{\infty\}$ in such a way that $\Delta_0$
corresponds to 
$\matH=\{(x,0)\in\matR^{n-1}\times\{0\}:\ x_{n-1}>0\}$.
Now let $p_1,p_2\in\Lambda(\Gamma)\setminus\partial \Delta_0$ and let 
$\alpha:[0,1]\to (\matR^{n-1}\times\{0\})\setminus\overline{\matH}$ be the straight Euclidean
segment which joins $p_1$ to $p_2$. If $\{(a_i,b_i)\subset[0,1],\ i\geqslant 1\}$
is the set of the connected components of $\alpha^{-1}(\Omega(\Gamma))$, 
then, up to reordering the $\Delta_i$'s with $i\geqslant 1$, we have
$\alpha([a_i,b_i])\subset\overline{\Delta}_i$. Let $r_i$ be the Euclidean
radius of $\Delta_i$. Since $\partial \Delta_i$ can touch $\partial \Delta_0$
at most in one point, for any $i\geqslant 1$ there exists a path 
$\beta_i:[a_i,b_i]\to\partial \Delta_i$ with $\beta_i(a_i)=\alpha(a_i)$,
$\beta_i(b_i)=\alpha(b_i)$ and $\mathrm{length}(\beta_i)\leqslant
2\pi r_i$. Now let $\alpha_i$ be the path inductively defined as follows:
$\alpha_0=\alpha$, $\alpha_{i+1}(t)=\beta_{i+1}(t)$ if $t\in [a_{i+1},b_{i+1}]$
and $\alpha_{i+1}(t)=\alpha_i (t)$ if $t\in [0,a_{i+1}]\cup [b_{i+1},1]$.
The path $\alpha_i$ is obviously continuous for any $i\in\matN$. Moreover,
since $\lim_{i\to\infty} r_i=0$, the sequence of paths $\{\alpha_i,i\in\matN\}$
uniformly converges to the desired continuous path 
$\alpha_{\infty}:[0,1] \to \Lambda(\Gamma)\setminus\partial \Delta_0$.

Suppose now that $\Lambda(\Gamma)\setminus j(S^{n-2})$ is path connected.
The Jordan-Brower separation theorem implies that $\partial \matH^n\setminus
j(S^{n-2})=A_1\cup A_2$, where the $A_i$'s are disjoint open subset 
of $\partial \matH^n$ 
with $\partial A_i=j(S^{n-2})$ for $i=1,2$ (since we are not assuming that
$j$ is tame, at this stage we are not allowed
to claim that the $A_i$'s are topological
balls). Our hypothesis now
forces $A_k\cap\Lambda(\Gamma)=\emptyset$ for some $k\in\{1,2\}$, so
$A_k\subset\Delta_l$ for some $l\in\matN$. Moreover, since
$\partial A_k=j(S^{n-2})\subset\Lambda(\Gamma)$, it is easily seen
that $j(S^{n-2})=\partial \Delta_l$, and we are done.
\finedimo

Form now on let $N_1$ and $N_2$ be hyperbolic $n$-manifolds, let
$\pi_i:\matH^n\supset\widetilde{N}_i\to N_i$ be the universal
covering of $N_i$ and let
$\Gamma_i$ be a discrete  
subgroup of $\mathrm{Iso}^+(\matH^n)$ such that 
$N_i\cong\widetilde{N}_i/\Gamma_i$.
Let also $\varphi:\Gamma_1\to\Gamma_2$ be a group 
isomorphism satisfying the condition of Theorem~\ref{main:teo}.
If $f:N_1\to N_2$ is a continuous map, it is easily seen
that $\varphi$ is induced by $f$ if and only if $f$
admits a continuous lift $\widetilde{f}:\widetilde{N}_1\to
\widetilde{N}_2$ such that $\widetilde{f}\circ\gamma=
\varphi(\gamma)\circ\widetilde{f}$ for every 
$\gamma\in\Gamma_1$.

\begin{lemma}\label{limitomeo:lem}
There exists a homeomorphism $\widehat{\varphi}:\Lambda(\Gamma_1)
\to\Lambda(\Gamma_2)$ such that 
$\widehat{\varphi}(\gamma(x))=\varphi(\gamma)(\widehat{\varphi}(x))$ for any
$x\in\Lambda(\Gamma_1)$, $\gamma\in\Gamma_1$.
\end{lemma}
\noindent\emph{Proof:}
For any group $G$, let us denote by $\overline{G}$ 
the \emph{completion} of $G$ (see~\cite{Floyd}
for a definition). Recall that $G$ acts in a natural
way on $\overline{G}$ as a group of homeomorphsims.
It is proved in~\cite{Floyd}
that any group isomorphism $\psi:G_1 \to G_2$
induces a homeomorphism $\overline{\psi}:
\overline{G}_1\to\overline{G}_2$ such that
$\overline{\psi}(g(x))=\psi(g)(\overline{\psi}(x))$.
Moreover, if $G$ is a geometrically finite subgroup
of $\mathrm{Iso}^+(\matH^n)$ then there exists a natural
continuous surjection $p_G:\overline{G}\to\Lambda(G)$ which
is 2-to-1 onto points with rank-1 parabolic stabilizer,
and injective everywhere else (this was shown in~\cite{Floyd} 
under the assumption $n=3$, but as it was observed in~\cite{Tukia}
the proof in~\cite{Floyd} actually works in any dimension).

Now  $\varphi$
induces by hypothesis a bijective correspondence
between rank-1 maximal parabolic subgroups of $\Gamma_1$ 
and rank-1 maximal parabolic subgroups
of  $\Gamma_2$. Using this fact it is easily seen
that there exists a unique bijective map
$\widehat{\varphi}:\Lambda(\Gamma_1)\to\Lambda(\Gamma_2)$ such that
$\widehat{\varphi}\circ p_{\Gamma_1}=p_{\Gamma_2}\circ\overline{\varphi}$.
Since $\overline{\Gamma}_i$ and $\Lambda(\Gamma_i)$ are
Haussdorff compact spaces for $i=1,2$, the map $\widehat{\varphi}$
is a homeomorphism, and we are done.
\finedimo

\begin{cor}\label{nobound:cor}
$\partial N_1=\emptyset$ if and only if
$\partial N_2=\emptyset$.
\end{cor}
\noindent\emph{Proof:}
Since $\Gamma_i$ is geometrically finite,
the boundary of $N_i$ is empty if and only
if $\Lambda(\Gamma_i)$ is homeomorphic to 
$S^{n-1}$. Lemma~\ref{limitomeo:lem} provides a 
homeomorphism between $\Lambda(\Gamma_1)$
and $\Lambda(\Gamma_2)$, and the conclusion
follows at once.
\finedimo

If $\partial N_1=\partial N_2=\emptyset$, 
Mostow-Prasad's rigidity theorem applies
ensuring geometricity of $\varphi$.
Then from now on we shall assume
that both $N_1$ and $N_2$ have non-empty 
boundary. 

\begin{lemma}\label{condiz:lem}
The isomorphism $\varphi$ satisfies the following
conditions:
\begin{enumerate}
\item
$\varphi(H)$ is a peripheral subgroup of $\Gamma_2$ if and only
if $H$ is a peripheral subgroup of $\Gamma_1$; if so
we also have $\widehat{\varphi}(\Lambda(H))=\Lambda(\varphi(H))$;
\item
$\varphi(\gamma)$ is a parabolic element of $\Gamma_2$
if and only if $\gamma$ is a parabolic element of
$\Gamma_1$.
\end{enumerate}
\end{lemma}
\noindent\emph{Proof:}
Let $\widehat{\varphi}:\Lambda(\Gamma_1)\to\Lambda(\Gamma_2)$ be the
homeomorphism constructed in Lemma~\ref{limitomeo:lem} and
let $H=\mathrm{stab}(\Delta)$ be a peripheral subgroup of 
$\Gamma_1$, where $\Delta$ is a component of $\Omega(\Gamma_1)$. 
By Lemma~\ref{connect:lemma}, $\Lambda(\Gamma_1)\setminus
\Lambda(H)=\Lambda(\Gamma_1)\setminus
\partial\Delta$ is path connected, so $\Lambda(\Gamma_2)
\setminus \widehat{\varphi}(\Lambda(H))=\widehat{\varphi}
(\Lambda(\Gamma_1)\setminus\Lambda(H))$ is also path connected,
and $\widehat{\varphi}(\Lambda(H))$ is equal to $\Lambda(K)$ for
some peripheral subgroup $K$ of $\Gamma_2$. Let $K=\mathrm{stab}
(\Delta')$, where $\Delta'$ is a component of $\Omega(\Gamma_2)$.
Now let $h$ be a loxodromic element of $H$ with fixed points
$p_1,p_2$ in $\Lambda(H)$. Since $\widehat{\varphi}$ is $\varphi$-equivariant,
we have that $\varphi(h)$ is a loxodromic element of $\Gamma_2$ 
with fixed points $\widehat{\varphi}(p_1), 
\widehat{\varphi}(p_2)$ which lie in $\Lambda(K)$.
Since the boundaries of two
different components of $\Omega(\Gamma_2)$ can intersect at most in one 
point, it easily follows that $\varphi(h)\in\mathrm{stab}(\Delta')=K$.
Now $H$ is generated by its loxodromic elements, so
$\varphi(H)$ is contained in $K$. On the other hand,
the same argument applied to $\varphi^{-1}$ shows that $\varphi^{-1}(K)$
is contained in a peripheral subgroup of $\Gamma_1$, say $H'$,
with $H\subset H'$. Now 
Remark~\ref{nomaximal:rem} implies that $H=H'$, so  $\varphi(H)=K$ 
and point~(1) is proved.

To prove point~(2), we observe that the $\varphi$-equivariance
of $\widehat{\varphi}$ implies that for any $\gamma\in\Gamma_1$ 
the fixed points of $\varphi(\gamma)$ are exactly
the images under $\widehat{\varphi}$ of the fixed points of $\gamma$.
This implies that the number of fixed points of $\varphi(\gamma)$
on $\Lambda(\Gamma_2)$ equals the number of fixed points
of $\gamma$ on $\Lambda(\Gamma_1)$, so $\varphi(\gamma)$
is parabolic if and only if $\gamma$ is.
\finedimo

\section{The $n$-dimensional case, $n\geqslant 4$}

The next proposition easily implies Theorem~\ref{main:teo}
under the assumption that the dimension of $N_1$ and $N_2$ is at least $4$.

\begin{prop}\label{final:prop}
Let $n\geqslant 4$. Then there exists a conformal map $f:\partial \matH^n
\to \partial \matH^n$ such that $f\circ\gamma=\varphi(\gamma)\circ
f$ for any $\gamma\in\Gamma_1$.
\end{prop}
\noindent\emph{Proof:}
Let $\Delta^{\!1}$ be a connected component of
$\Omega(\Gamma_1)$, and $H_1$ be the
stabilizer of $\Delta^{\!1}$ in $\Gamma_1$.
By Lemma~\ref{condiz:lem}, the group $H_2=\varphi(H_1)$ 
is a peripheral subgroup of $\Gamma_2$. 
Let now $\Delta^{\!2}$ be 
the $H_2$-invariant component of $\Omega(\Gamma_2)$,
\emph{i.e.} the unique component of $\Omega(\Gamma_2)$ whose boundary
is equal to $\Lambda(H_2)$.  
By Lemma~\ref{condiz:lem}~(1), the homeomorphism constructed in 
Lemma~\ref{limitomeo:lem} restricts to a homeomorphism
$\widehat{\varphi}|_{\partial \Delta^{\!1}}:\partial\Delta^{\!1}
\to\partial \Delta^{\!2}$
such that $\widehat{\varphi}|_{\partial \Delta^{\!1}}\circ\gamma
=\varphi(\gamma)\circ\widehat{\varphi}|_{\partial \Delta^{\!1}}$
for every $\gamma\in H_1$. Let now $S^1,S^2$ be the hyperplanes of $\matH^n$
bounded respectively by $\partial\Delta^{\!1}$ and $\partial\Delta^{\!2}$. Then
$S^k/H_k$ is isometric to a component of the geodesic
boundary of $N_k$ for $k=1,2$, so it is a finite-volume complete hyperbolic
$(n-1)$-manifold without boundary. Since $n\geqslant 4$, Mostow-Prasad's
rigidity theorem applies providing an isometry $g:S^1\to S^2$
whose continuous extension to $\partial \Delta^{\!1}$ is equal 
to $\widehat{\varphi}|_{\partial \Delta^{\!1}}$. 
%and such that $g\circ\gamma=\varphi(\gamma)\circ g$ for any $\gamma\in
%\Gamma_1^0$. 
Let now $p_k,k=1,2$ be the orthogonal projection of
$S^k$ onto $\Delta^{\!k}$, \emph{i.e.} the function which maps a point
$q\in S^k$ to the point $p\in\Delta^{\!k}$ such that the geodesic ray $[q,p)$ is orthogonal  
to $S^k$. The map $g':\Delta^{\!1}\to \Delta^{\!2}$
defined by $g'=p_2\circ g\circ p_1^{-1}$ is conformal, and
its continuous extension to $\partial \Delta^{\!1}$ is equal to $\widehat{\varphi}|_{\partial \Delta^{\!1}}$. 
%verifies the same $\varphi$-equivariance property as $g$.

By repeating the construction described above for each component of $\Omega(\Gamma_1)$, 
we can construct a conformal map $t:\Omega(\Gamma_1)\to\Omega(\Gamma_2)$. 
This map is a homeomorphism, since it admits
a continuous inverse which can be constructed from the isomorphism
$\varphi^{-1}:\Gamma_2\to \Gamma_1$. We want now to show that for any
$\gamma\in\Gamma_1$, we have $t\circ\gamma=\varphi(\gamma)\circ t$.
Let $\Delta$ be a component of $\Omega(\Gamma_1)$.
By the very definition of $t$ it follows that $t(\Delta)$ is the unique component
of $\Omega(\Gamma_2)$ which is bounded by $\widehat{\varphi}(\partial \Delta)$, so
\begin{eqnarray*}
\partial (\varphi(\gamma)(t(\Delta)))&=&\varphi(\gamma)(\partial (t(\Delta)))=
\varphi(\gamma)(\widehat{\varphi}(\partial \Delta))=\widehat{\varphi}(\gamma(\partial \Delta))\\
&=&
\widehat{\varphi}(\partial(\gamma(\Delta)))=
\partial (t(\gamma(\Delta)).
\end{eqnarray*} 
This shows that both $t\circ\gamma$ and $\varphi(\gamma)\circ t$ map $\Delta$ onto
the same component $\Delta'$ of $\Omega(\Gamma_2)$. Moreover, the continuous extensions
of $t\circ\gamma$ and $\varphi(\gamma)\circ t$ to $\partial \Delta$ 
are respectively equal to $\widehat{\varphi}\circ\gamma$ and 
$\varphi(\gamma)\circ\widehat{\varphi}$, which are in turn equal to each other 
because of the $\varphi$-equivariance of $\widehat{\varphi}$.
Being conformal, the maps $t\circ\gamma$ and $\varphi(\gamma)\circ t$ 
must then be equal on $\Delta$, and this proves 
the required $\varphi$-equivariance of $t$.

Now let $f:\partial \matH^n\to\partial\matH^n$ be defined by
$f(x)=t(x)$ if $x\in\Omega(\Gamma_1)$, and $f(x)=\widehat{\varphi}(x)$ 
if $x\in\Lambda(\Gamma_1)$. To conclude the proof we only have to 
observe that since $f$ is $\varphi$-equivariant and conformal on
$\Omega(\Gamma_1)$, a result of Tukia~\cite{Tukia2}
ensures that $f$ is a coformal map.
\finedimo

We can now conclude the proof of Theorem~\ref{main:teo}, under the
assumption that the dimension of $N_1$ and $N_2$ is greater than 3. 
Let $\widetilde{\psi}$ be the unique isometry of $\matH^n$ whose
continuous extension to $\partial \matH^n$ is equal to $f$.
The $\varphi$-equivariance of $f$ readily implies that
$\widetilde{\psi}(\gamma(x))=\varphi(\gamma)(\widetilde{\psi}(x))$ for every 
$x\in\matH^n,\gamma\in\Gamma_1$. If we identify $N_i$ with the convex core
of the manifold $\matH^n/\Gamma_i$ for $i=1,2$, then $\widetilde{\psi}$ induces an isometry  
$\psi:N_1\to N_2$ with $\psi_{\ast}=\varphi$.

\section{The 3-dimensional case}

As briefly explained in the introduction, the 3-dimensional case needs
a different approach.  

\begin{lemma}\label{final:lem}
There exists a homeomorphism $g:N_1\to N_2$ such that
$\varphi=g_\ast$.
\end{lemma}
\noindent\emph{Proof:}
Let $M_i=\left(\matH^3\cup\Omega(\Gamma_i)\right)/\Gamma_i$ for $i=1,2$.
By Lemma~\ref{condiz:lem} and Remark~\ref{nomaximal:rem}, we can apply
Theorem~1 of~\cite{MM} to $\varphi$, obtaining a homeomorphism
$g':M_1\to M_2$ inducing $\varphi$ (note that our definition
of \emph{geometric} is stronger than the
one in~\cite{MM}). Now $N_i$ is canonically embedded
in $M_i$ in such a way that $M_i\setminus N_i$ is an open
collar of $\partial M_i$. This implies that $g'$ can be isotoped
to a $g'':M_1\to M_2$ such that $g''(N_1)=N_2$
and $g=g''|_{N_1}$ is the required homeomorphism.
\finedimo

\begin{rem}\label{joha:rem}
\emph{If $N_1$ and $N_2$ have compact geodesic boundary, then Lemma~\ref{final:lem}
can also be deduced by the following result of Johannson~\cite{joha,swarup}:
Any homotopy equivalence between compact orientable 
boundary-irreducible anannular Haken
$3$-manifolds can be homotoped to a homeomorphism.}
\end{rem}
 
We can now conclude the proof of Theorem~\ref{main:teo} in the case
when $N_1$ and $N_2$ are $3$-dimensional manifolds.
Let $g:N_1 \to N_2$ be the homeomorphism constructed in
Lemma~\ref{final:lem}, let $D(N_i)$ be the double of $N_i$ for
$i=1,2$ and let $D(g):D(N_1)\to D(N_2)$ be the homeomorphism
obtained by doubling $g$. By Mostow-Prasad's rigidity theorem,
$D(g)$ is homotopic to an isometry $h:D(N_1)\to D(N_2)$.
Since $\partial N_2=g(\partial N_1)$ and $h(\partial N_1)$ are embedded
totally geodesic homotopic surfaces in $N_2$, 
we get that $h(\partial N_1)=\partial N_2$, so $h(N_1)=N_2$.
Moreover, $h_{\ast}=g_{\ast}$ on $\pi_1(D(N_1))$, and 
the inclusion of $\pi_1(N_i)$ in $\pi_1(D(N_i))$ is injective
for $i=1,2$, so $h_{\ast}=g_{\ast}=\varphi$
on $\Gamma_1$. In conclusion, we have shown that $h|_{N_1}:
N_1\to N_2$ is an isometry inducing $\varphi$, so $\varphi$
is geometric.

\paragraph{Counterexamples in the non-compact boundary case}
We now show that the conclusions of Corollaries~\ref{cor1}
and~\ref{cor2}
are no longer true if we consider hyperbolic
$3$-manifolds with non-compact geodesic boundary.
More precisely, we will prove the following:

\begin{prop}\label{prop:contro}
There exist hyperbolic $3$-manifolds with non-compact geodesic boundary
$N_1,N_2$ such that:
\begin{enumerate}
\item
$\pi_1(N_1)\cong\pi_1(N_2)$ but $N_1$ is not homeomorphic to $N_2$;
\item
$\mathrm{Out}(\pi_1(N_i))\ncong\mathrm{Iso}(N_i)$ for $i=1,2$.
\end{enumerate}
\end{prop}
\noindent\emph{Proof:}
We will give an explicit construction of $N_1$ and $N_2$.
\begin{figure}[ht]
\begin{center}
\input{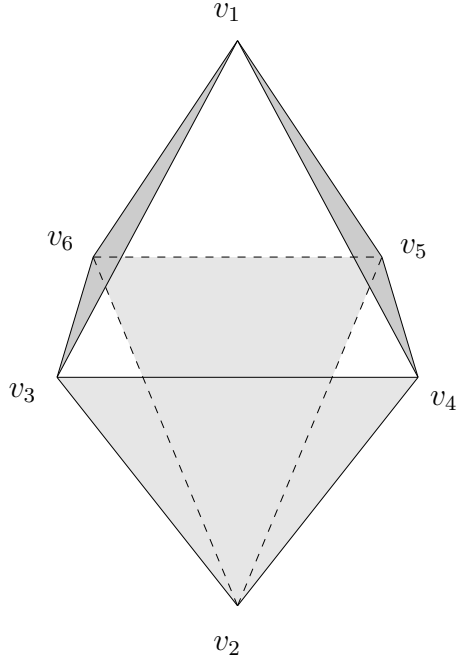}
\caption{\small{The manifolds $N_1,N_2$ and $N_3$ 
are obtained by gluing in pairs the non-shadowed
faces of the regular ideal octahedron 
along suitable isometries.}}\label{octa:fig}
\end{center}\end{figure}
Let $O\subset\matH^3$ be the regular ideal octahedron and let
$v_1,\ldots,v_6$ be the vertices of $O$ as shown in Fig.~\ref{octa:fig}.
We denote by $F_{ijk}$ the face of $O$ with vertices 
$v_i,v_j,v_k$. Let $g:F_{134}\to F_{156}$ be the unique 
orientation-reversing isometry such that $g(v_1)=v_6$, and 
$h_1,h_2:F_{542}\to F_{362}$ be the unique orientation-reversing
isometries such that $h_1(v_5)=v_6$, $h_2(v_5)=v_3$.
We now define $N_1$ to be the manifold obtained by gluing $O$
along $g$ and $f_1$, and $N_2$ to be the manifold 
obtained by gluing $O$ along $g$ and $f_2$.   
Since all the dihedral angles of $O$ are right, it is 
easily seen that the metric on $O$ induces a complete
finite-volume hyperbolic structure on the $N_i$'s 
such that the shadowed 
faces in Fig.~\ref{octa:fig} are glued along their egdes to give
a non-compact totally geodesic boundary.

Now the natural compactification of $N_i$ is homeomorphic 
to the genus-2 handlebody for $i=1,2$, so
$\pi_1(N_1)\cong\pi_1(N_2)\cong \matZ\ast\matZ$. 
Moreover, the boundary of
$N_1$ is homeomorphic to the $2$-punctured torus, 
while the boundary of $N_2$ is homeomorphic to the 
$4$-punctured sphere, so $N_1$ is not homeomorphic to 
$N_2$. This proves point~(1).

In order to prove point~(2), we only have to observe
that the group of the outer isomorphisms of $\matZ\ast\matZ$ is of infinite
order, while the group of isometries of any complete
finite-volume hyperbolic
$n$-manifold with geodesic boundary has a finite
number of elements.
\finedimo

\begin{ex}\label{omeobordo:rem}
\emph{Let $N_3$ be the hyperbolic manifold with non-compact geodesic boundary
obtained by gluing the faces of $O$ along $h_2$ and $g'$, where 
$g':F_{134}\to F_{156}$ is the unique orientation-reversing isometry
such that $g'(v_1)=v_5$. As before, the natural
compactification of $N_3$ is the genus-$2$ handlebody,
so $\pi(N_3)\cong\pi(N_2)\cong \matZ\ast\matZ$.
Moreover, with some effort one could show that $\partial N_2$ is 
homeomorphic but not isometric to $\partial N_3$,
and $N_2$ and $N_3$ 
are not homeomorphic to each other.} 
\end{ex}

\paragraph{A more general construction}
We now briefly describe a different method of contructing
homotopically-equivalent non-homeomorphic 
hyperbolic $3$-manifolds with non-compact geodesic boundary.
To this aim we first 
recall that Thurston's hyperbolization theorem
for Haken manifolds~\cite{thu2} gives necessary and sufficient topological
conditions on a manifold to be hyperbolic with
geodesic boundary:

\begin{teo}\label{hyper:teo}
Let $\overline{M}$ be a compact orientable manifold 
with non-empty boundary, let $\calT$ be the set of
boundary tori of $\overline{M}$ and
let $\calA$ be a family of disjoint closed annuli in 
$\partial \overline{M}\setminus\calT$. 
Then $M=\overline{M}\setminus (\calT\cup\calA)$ is hyperbolic
if and only if the pair $(\overline{M},\calA)$ satisfies
the following conditions:
\begin{itemize}
\item
the components of $\partial M$ 
have negative Euler characteristic;
\item
$\overline{M}\setminus\calA$ 
is boundary-irreducible and geometrically atoroidal;
\item
the only proper essential annuli contained in $M$ are parallel
in $\overline{M}$ to the annuli in $\calA$.
\end{itemize}
\end{teo}

Using Theorem~\ref{hyper:teo} we will now prove
the following:

\begin{prop}\label{contro:prop}
Let $N$ be a hyperbolic $3$-manifold with non-empty
geodesic boundary, and suppose that at least one component
of $\partial N$ is not a $3$-punctured sphere.
Then there exists a hyperbolic $3$-manifold with geodesic boundary
which is homotopically equivalent but not homeomorphic to $N$.
\end{prop}
\noindent\emph{Proof:}
Let $\overline{N}$ be the natural compactification of $N$ 
obtained by adding to $N$ a family 
$\calA_N$ of closed annuli and a family
$\calT_N$ of tori.
Let also $\{\alpha_1,
\ldots,\alpha_k\}$ be a non-empty 
family of disjoint essential non-parallel
loops on $\partial N$ (such a family always exists
because of the assumption on $\partial N$). 
Let $\calA'$ be the family of annuli
in $\partial \overline{N}\setminus\calT_N$ obtained by adding to
$\calA_N$ closed regular neighbourhoods in $\partial N$ of the
$\alpha_i$'s. It is easily seen that the pair $(\overline{N},\calA')$
satisfies the conditions of Theorem~\ref{hyper:teo},
so $N'=N\setminus(\bigcup_{i=1}^k \alpha_i)$ is hyperbolic.
Of course $N'$ is homotopically equivalent to $N$, but
$\partial N'$ is not homeomorphic
to $\partial N$, so \emph{a fortiori} $N$ and $N'$ are 
not homeomorphic to each other.
\finedimo

\vspace{1.5 cm}

\noindent
\hspace*{6cm}Scuola Normale Superiore\\ 
\hspace*{6cm}Piazza dei Cavalieri 7 \\
\hspace*{6cm}56127 Pisa, Italy\\ 
\hspace*{6cm}frigerio@sns.it
\vspace{.5 cm}

%\noindent
%\hspace*{6cm}Dipartimento di Matematica\\ 
%\hspace*{6cm}Universit\`a di Pisa\\ 
%\hspace*{6cm}Via F. Buonarroti 2\\ 
%\hspace*{6cm}56127 Pisa, Italy\\ 
%\hspace*{6cm}martelli@mail.dm.unipi.it
%\vspace{.5 cm} 

%\noindent
%\hspace*{6cm}Dipartimento di Matematica Applicata\\ 
%\hspace*{6cm}Universit\`a di Pisa\\
%\hspace*{6cm}Via Bonanno Pisano 25B\\ 
%\hspace*{6cm}56126 Pisa, Italy\\ 
%\hspace*{6cm}petronio@dm.unipi.it


\begin{thebibliography}{99}

\bibitem{Floyd} \textsc{W.J.~Floyd},
\emph{Group completions and limit sets of Kleinian groups},
Invent. Math. \textbf{57} (1980), 205-218.

\bibitem{FriPe} \textsc{R.~Frigerio, C.~Petronio},
\emph{Construction and recognition of hyperbolic
manifolds with geodesic boundary},
{\tt math.GT/0109012}, to appear in Trans. Amer. Math. Soc.

\bibitem{joha} \textsc{K.~Johannson},
\emph{Homotopy equivalences of $3$-manifolds with boundaries},
Lecture Notes in Mathematics, 761. Springer, Berlin, 1979. 


\bibitem{KMS} \textsc{L.~Keen, B.~Maskit, C.~Series},
\emph{Geometric finiteness and uniqueness 
for Kleinian groups with circle packing limit sets}, 
J. Reine Angew. Math. \textbf{436} (1993), 209-219.

\bibitem {Kojima1} \textsc{S.~Kojima},
\emph{Polyhedral decomposition of hyperbolic $3$-manifolds
with totally geodesic boundary},
``Aspects of low-dimensional manifolds, Kinokuniya, Tokyo'',
Adv. Stud. Pure Math. \textbf{20} (1992), 93-112.

\bibitem{Kojima2} \textsc{S.~Kojima},
\emph{Geometry of hyperbolic $3$-manifolds with boundary},
Kodai Math. J. \textbf{17} (1994), 530-537.


%\bibitem{Marden} \textsc{A.~Marden},
%\emph{The geometry of finitely generated kleinian groups},
%Ann. of Math. (2) \textbf{99} (1974), 383-462.

\bibitem{MM} \textsc{A.~Marden, B.~Maskit}, 
\emph{On the isomorphism theorem for Kleinian groups},
Invent. Math. \textbf{51} (1979), 9-14. 

%\bibitem{KM} \textsc{Y.~Miyamoto},
%\emph{Volumes of hyperbolic manifolds with geodesic boundary},
%Topology \textbf{33} (1994), 613-629. 




\bibitem{rat} \textsc{J.~Ratcliffe},
\emph{Foundations of hyperbolic manifolds}, 
Graduate Texts in Mathematics, 149. Springer-Verlag, New York, 1994.

\bibitem{swarup} \textsc{G.A.~Swarup},  
\emph{On a theorem of Johannson}, 
J. London Math. Soc. (2) \textbf{18} (1978), 560-562.


\bibitem{thu} \textsc{W.P.~Thurston},
``The geometry and topology of $3$-manifolds'',
mimeographed notes, Princeton, 1979. 


\bibitem{thu2} \textsc{W.P.~Thurston},
\emph{Three-dimensional manifolds, Kleinian groups and hyperbolic geometry},
Bull. Amer. Math. Soc. (N.S.) \textbf{6} (1982), 357-381.

\bibitem{Tukia2} \textsc{P.~Tukia},
\emph{On isomorphisms of geometrically finite Moebius groups},
Inst. Hautes \'Etudes Sci. Publ. Math. \textbf{61} (1985), 171-214. 

%\bibitem{Tukia3} \textsc{P.~Tukia,~J.~V\"ais\"al\"a},
%\emph{A remark on $1$-quasiconformal maps}, 
%Ann. Acad. Sci. Fenn. Ser. A I Math. \textbf{10} (1985), 561-562. 

\bibitem{Tukia} \textsc{P.~Tukia},
\emph{A remark on a paper by Floyd}, in 
``Holomorphic functions and moduli, Vol. II'' (Berkeley, CA, 1986), 165--172, 
Math. Sci. Res. Inst. Publ., \textbf{11}, Springer, New York, 1988.

%\bibitem{Wald} \textsc{F.~Waldhausen}, 
%\emph{On irreducible 3-manifolds which are sufficiently large},
%Ann. of Math. (2) \textbf{87} (1968), 56-88.



\end{thebibliography}
\end{document}